\documentclass[12pt]{article}
\usepackage{amssymb}
\usepackage{amsfonts}
\usepackage{latexsym}
\usepackage{amsthm}

\theoremstyle{plain}
\newtheorem{theorem}{\bf Theorem}[section]
\newtheorem{lemma}[theorem]{\bf Lemma}

\theoremstyle{remark}

\newcommand{\er}[1]{{\rm(\ref{#1})}}
\def\lb{\label}

\begin{document}
\def\a{\alpha}  \def\cA{{\cal A}}    \def\bA{{\bf A}}  \def\mA{{\mathscr A}}
\def\b{\beta}   \def\cB{{\cal B}}    \def\bB{{\bf B}}  \def\mB{{\mathscr B}}
\def\g{\gamma}  \def\cC{{\cal C}}    \def\bC{{\bf C}}  \def\mC{{\mathscr C}}
\def\G{\Gamma}  \def\cD{{\cal D}}    \def\bD{{\bf D}}  \def\mD{{\mathscr D}}
\def\d{\delta}  \def\cE{{\cal E}}    \def\bE{{\bf E}}  \def\mE{{\mathscr E}}
\def\D{\Delta}  \def\cF{{\cal F}}    \def\bF{{\bf F}}  \def\mF{{\mathscr F}}
\def\c{\chi}    \def\cG{{\cal G}}    \def\bG{{\bf G}}  \def\mG{{\mathscr G}}
\def\z{\zeta}   \def\cH{{\cal H}}    \def\bH{{\bf H}}  \def\mH{{\mathscr H}}
\def\e{\eta}    \def\cI{{\cal I}}    \def\bI{{\bf I}}  \def\mI{{\mathscr I}}
\def\p{\psi}    \def\cJ{{\cal J}}    \def\bJ{{\bf J}}  \def\mJ{{\mathscr J}}
\def\vT{\Theta} \def\cK{{\cal K}}    \def\bK{{\bf K}}  \def\mK{{\mathscr K}}
\def\k{\kappa}  \def\cL{{\cal L}}    \def\bL{{\bf L}}  \def\mL{{\mathscr L}}
\def\l{\lambda} \def\cM{{\cal M}}    \def\bM{{\bf M}}  \def\mM{{\mathscr M}}
\def\L{\Lambda} \def\cN{{\cal N}}    \def\bN{{\bf N}}  \def\mN{{\mathscr N}}
\def\m{\mu}     \def\cO{{\cal O}}    \def\bO{{\bf O}}  \def\mO{{\mathscr O}}
\def\n{\nu}     \def\cP{{\cal P}}    \def\bP{{\bf P}}  \def\mP{{\mathscr P}}
\def\r{\rho}    \def\cQ{{\cal Q}}    \def\bQ{{\bf Q}}  \def\mQ{{\mathscr Q}}
\def\s{\sigma}  \def\cR{{\cal R}}    \def\bR{{\bf R}}  \def\mR{{\mathscr R}}
\def\S{\Sigma}  \def\cS{{\cal S}}    \def\bS{{\bf S}}  \def\mS{{\mathscr S}}
\def\t{\tau}    \def\cT{{\cal T}}    \def\bT{{\bf T}}  \def\mT{{\mathscr T}}
\def\f{\phi}    \def\cU{{\cal U}}    \def\bU{{\bf U}}  \def\mU{{\mathscr U}}
\def\F{\Phi}    \def\cV{{\cal V}}    \def\bV{{\bf V}}  \def\mV{{\mathscr V}}
\def\P{\Psi}    \def\cW{{\cal W}}    \def\bW{{\bf W}}  \def\mW{{\mathscr W}}
\def\o{\omega}  \def\cX{{\cal X}}    \def\bX{{\bf X}}  \def\mX{{\mathscr X}}
\def\x{\xi}     \def\cY{{\cal Y}}    \def\bY{{\bf Y}}  \def\mY{{\mathscr Y}}
\def\X{\Xi}     \def\cZ{{\cal Z}}    \def\bZ{{\bf Z}}  \def\mZ{{\mathscr Z}}
\def\O{\Omega}
\def\ve{\varepsilon}
\def\eps{\epsilon}
\def\vt{\vartheta}
\def\vp{\varphi}
\def\vk{\varkappa}

\def\Z{{\Bbb Z}}
\def\R{{\Bbb R}}
\def\C{{\Bbb C}}
\def\T{{\Bbb T}}
\def\N{{\Bbb N}}

\let\ge\geqslant
\let\le\leqslant
\def\ma{\left(\begin{array}{cc}}
\def\am{\end{array}\right)}
\def\iint{\int\!\!\!\int}
\def\lt{\biggl}
\def\rt{\biggr}
\let\geq\geqslant
\let\leq\leqslant
\def\[{\begin{equation}}
\def\]{\end{equation}}
\def\wt{\widetilde}
\def\pa{\partial}
\def\sm{\setminus}
\def\es{\emptyset}
\def\no{\noindent}
\def\ol{\overline}
\def\iy{\infty}
\def\ev{\equiv}
\def\/{\over}
\def\ts{\times}
\def\os{\oplus}
\def\ss{\subset}
\def\h{\hat}
\def\Re{\mathop{\rm Re}\nolimits}
\def\Im{\mathop{\rm Im}\nolimits}
\def\supp{\mathop{\rm supp}\nolimits}
\def\sign{\mathop{\rm sign}\nolimits}
\def\Ran{\mathop{\rm Ran}\nolimits}
\def\Ker{\mathop{\rm Ker}\nolimits}
\def\Tr{\mathop{\rm Tr}\nolimits}
\def\const{\mathop{\rm const}\nolimits}
\def\Wr{\mathop{\rm Wr}\nolimits}
\def\Ai{\mathop{\rm Ai}\nolimits}
\def\Bi{\mathop{\rm Bi}\nolimits}
\def\BBox{\hspace{1mm}\vrule height6pt width5.5pt depth0pt \hspace{6pt}}

\def\Twelve{
\font\Tenmsa=msam10 scaled 1200 \font\Sevenmsa=msam7 scaled 1200
\font\Fivemsa=msam5 scaled 1200
%\newfam\msafam
\textfont\msafam=\Tenmsa \scriptfont\msafam=\Sevenmsa
\scriptscriptfont\msafam=\Fivemsa \font\Tenmsb=msbm10 scaled 1200
\font\Sevenmsb=msbm7 scaled 1200 \font\Fivemsb=msbm5 scaled 1200
%\newfam\msafam
\textfont\msbfam=\Tenmsb \scriptfont\msbfam=\Sevenmsb
\scriptscriptfont\msbfam=\Fivemsb

\font\Teneufm=eufm10 scaled 1200 \font\Seveneufm=eufm7 scaled 1200
\font\Fiveeufm=eufm5 scaled 1200
%\newfam\eufmfam
\textfont\eufmfam=\Teneufm \scriptfont\eufmfam=\Seveneufm
\scriptscriptfont\eufmfam=\Fiveeufm}

\def\Ten{
\textfont\msafam=\tenmsa \scriptfont\msafam=\sevenmsa
\scriptscriptfont\msafam=\fivemsa

\textfont\msbfam=\tenmsb \scriptfont\msbfam=\sevenmsb
\scriptscriptfont\msbfam=\fivemsb

\textfont\eufmfam=\teneufm \scriptfont\eufmfam=\seveneufm
\scriptscriptfont\eufmfam=\fiveeufm}

\title {Inverse problem for the discrete 1D Schr\"odinger operator with small periodic potentials}

\author{ Evgeny Korotyaev
\begin{footnote}
{Institut f\"ur  Mathematik,  Humboldt Universit\"at zu Berlin,
Rudower Chaussee 25, 12489, Berlin, Germany, e-mail:
evgeny@math.hu-berlin.de \ \ To whom correspondence should be
addressed}
\end{footnote}
and  Anton Kutsenko%${}^{\ast,\!\!\!}$
\begin{footnote}
{Faculty of Math. and Mech. St-Petersburg State University,\
%Institut f\"ur  Mathematik, 
%Universit\"at Potsdam, 
e-mail:kucenko@rambler.ru}
\end{footnote}
} \maketitle

\begin{abstract}
\no Consider the discrete 1D Schr\"odinger operator on $\Z$ with
an odd $2k$ periodic potential $q$. For small potentials we show
that the mapping: $q\to $ heights of vertical
slits on the quasi-momentum domain (similar to the
Marchenko-Ostrovski maping for the Hill operator) is a local
isomorphism and  the isospectral set consists of $2^k$ distinct
potentials. Finally, the asymptotics of the spectrum are
determined as $q\to 0$.
\end{abstract}

\vskip 0.25cm

\section {Introduction and main results}
\setcounter{equation}{0}

We consider the Schr\"odinger operator $(L
y)_n=y_{n-1}+y_{n+1}+q_ny_n,\ \ n\in\Z,$ acting on $l^2(\Z)$,
 where $\{q_n\}_{-\iy}^{\iy}$ is a real $N+1$ periodic sequence,
$q_{n+N+1}=q_n,\ n \in \Z$ and
\[
 \lb{1a} q\ev\{q_n\}_1^{N+1}\in \cQ\equiv\lt\{q\in\R^{N+1}:\
 \sum_1^{N+1}q_n=0\rt\},\ \ \|q\|^2=(q,q)=\sum_{n=1}^{N+1} q_n^2.
\]
It is well known that the spectrum of $L$ is absolutely continuous
and consists of $N+1$ intervals $\s_n=\s_n(q)=[\l_n^+,\l_{n+1}^-],
\ n=0,1, ... ,N$, where $\l_n^\pm =\l_n^\pm (q)$ and
$\l_0^{+}<\l_1^-\le \l_1^+< ... < \l_N^- \le \l_N^+ < \l_{N+1}^{-}
$. These intervals are separated by gaps
$\g_n=\g_n(q)=(\l_n^-,\l_n^+)$ of lengths $|\g_n|\ge 0$. If a gap
$\g_n$ is degenerate, i.e. $|\g_n|=0$, then the corresponding
segments    merge. Introduce fundamental solutions $\vp_n(\l,q)$
and $\vt_n(\l,q),\ {n\in \Z} $ of the equation
\[
 \lb{1b}
 y_{n-1}+y_{n+1}+q_ny_n=\l y_n,\ \ (\l,n)\in\C\ts\Z,
\]
with initial conditions $ \vp_{0}(\l,q)\ev \vt_1(\l,q)\ev 0,\ \ \
\vp_1(\l,q)\ev \vt_{0}(\l,q)\ev 1.$
 The function
$\D(\l,q)=\vp_{N+2}(\l,q)+\vt_{N+1}(\l,q)$ is called the Lyapunov
function for the operator $L$. The functions $\D, \vp_n$ and
$\vt_n$ are polynomials of $(\l,q)\in\C^{N+2}$.
 The spectrum of $L$ is given by
$\s(q)=\{\l\in\R:\ |\D(\l,q)|\le 2\}$ (see e.g. [Te]). Note that
$(-1)^{N+1-n}\D(\l_n^{\pm},q)=2, n=0,...,N+1$.  For each
$n=1,..,N$ there exists a unique $\l_n=\l_n(q)\in[\l_n^-,\l_n^+]$
such that
\[
\lb{1c}
 \D'(\l_n,q)=0,\ \ \ \D''(\l_n,q)\neq 0,\ \ \
 (-1)^{N+1-n}\D(\l_n,q)\ge2.
\]
Here and below we use the notation $(\ ')={\pa/\pa\l}$, ${\pa}_n=
{\pa/\pa q_n}$. Define $h_n, n=1,.., N$ by the equation
\[
\lb{1d}
 \D(\l_n,q)=2(-1)^{N+1-n}\cosh h_n,\ \ \  h_n\ge 0.
\]
The inverse spectral problem consists of the following parts:\\
i) Uniqueness.
Prove that the spectral data uniquely determine the potential.\\
ii) Characterization.
Give conditions for some data to be the spectral data of some potential.\\
iii) Reconstruction. Give an  algorithm for recovering the
potential from the spectral data.\\
 iv) A priori estimates (stability). Obtain estimates of the potential in terms of the  spectral data.

There is an enormous literature on the scalar Hill operator
including the inverse spectral theory. We mention the papers where
the inverse problem (including a characterization) was solved. In
the following papers  [MO1-2], [GT], [KK], [K1-2]  the authors
show that the mapping: $\{$potential$\}$ $\to $ $\{$spectral
data$\}$ is a homeomorphism. In particular this gives Uniqueness
and Characterization. In the recent paper [K3] one of the authors
solved the inverse problem (including a characterization) for the
case $-y''+v'y$, where $v\in L^2(\T)$ (i.e., $v'$ is a
distribution) and $\T=\R/\Z$. There are many papers about the
inverse problem for periodic Jacobi matrices (see ref. in
[Te]). But the corresponding extension to the case of periodic
Jacobi matrices was made only in [BGGK], [KKu1],[K4]. Remark that
the paper [vM] did not discuss the characterization of
spectral data for the periodic Jacobi matrices. However, in spite
of the importance of generalizing these studies in the continuous
case to the discrete 1D Schr\"odinger operator, until recently no
essential result (apart from the information given for the periodic
Jacobi matrices) has been proved. Note that for the 2D 
Schr\"odinger operator there is a book [GKT] about the inverse spectral problem.

The analysis of the discrete Schr\"odinger operator poses
interesting new problems: 1) to construct the mapping $q\to \bS$
(the spectral data) and solve the corresponding inverse problem
(we need a good choice of the spectral data $\bS$, a priori we
have a lot of candidates), 2) to study the quasimomentum (the real
part of the quasimomentum is the integrated density of states) as a
conformal mapping, 3) to obtain a priori estimates of potentials
in terms of spectral data, 4) recovering the potential using 
the spectral data.

The most complicated problem is Characterization, but even the
problem of uniqueness still is not solved, since we do not know
the "good" spectral data. 
Concerning Characterization, this paper gives a partial result 
for $q$ small (in the ball around the origin outside some cones of small volume ).
In fact this is the motivation of our
paper. Note that the case of large $q$ was studied in [KKu2], see
Theorem \ref{T3}.

Introduce the space of odd potentials $\cQ^{odd}$ by
\[
 \cQ^{odd}=\{q\in\R^{N+1}:\ q_{N+2-n}=-q_{n},\ n=1,...,N\},\ \
 N+1=2k,\ \ k\in\N.
\]
Define the mapping $ h:\cQ^{odd}\to \R^k$ by $q\to h(q)={\{
h_n(q)\}}_{1}^k,$ where $h_n$ is given by \er{1d}. The mapping $
h$ is some analog of the Marchenko-Ostrovski mapping for the Hill
operator [M], [MO1],[MO2]. The parameters $h_n$ have also a geometric sense as
the height of the vertical slits of the quasimomentum domain (see e.g.
[KKr]).  For given $ h(q)$ for some $q\in \cQ^{odd}$, the Lyapunov
function $\D(\l,q)$ is uniquely determined (see [Pe], [KKu1]).
Note that in our case of odd potentials in order to determine
$\D(\l,q)$ we only need heights $h_n,\ n=1,..,k$, since due to
\er{iL00} the  function $\D(\cdot,q),\ q\in\cQ^{odd}$ is even . If
$q,p\in \cQ^{odd}$, then
\[
\lb{pe1}
 \s(q)=\s(p)\ \ \ \Leftrightarrow\ \ \ \ \ \D(\cdot,q)=\D(\cdot,p)
\ \ \     \Leftrightarrow \ \ \  h(q)=h(p) .
\]
This simple fact follows from the result for Jacobi matrices (see
[Pe], [KKu1]). Introduce an orthonormal basis $\hat{e}_m$ and
coordinates $\hat q_m,m=1,..k$ in the space $\cQ^{odd}$ by
\[
\lb{17}
 \hat{e}_m=\frac{2^{-\frac{\d_{k,m}}2}}{\sqrt{k}}
 \left\{\sin \frac{(2n-1)m\pi}{2k}\right\}_{n=1}^{N+1},\ \
 (\hat{e}_m,\hat{e}_l)=\d_{m,l},\ \ m,l=1,...,k, \  q=\sum_{m=1}^k\hat q_m\hat e_m,\ \
\]
where $(\cdot,\cdot)$ is a standard  scalar product in $\C^{N+1}, N+1=2k$ and $\d_{k,m}$ is the Kroneker symbol.
Recall that $\l_n^{\pm}(0)=-2\cos\frac{\pi n}{N+1}, n=1,..,N$. We
formulate our first result.

\begin{theorem}  \lb{T1}
Let $q\in\cQ^{odd}$ and $\|q\|\to 0$. Then for any $n=1,...,N$ the
following asymptotics are fulfilled:
\[
 \lb{i0}
 (\l_n^{\pm}(q)-\l_n^{\pm}(0))^2=
 \frac{2^{\d_{k,n}}}{4k}\hat q_n^2+O(\|q\|^3),\ \ \ \
\]
\[
 \lb{9976}
 h_n^2(q)=\frac{2^{\d_{k,n}}k}{4\sin^2\frac{\pi n}{2k}}
 \hat q_n^2+O(\|q\|^4).
\]
\end{theorem}

Recall that the asymptotics for $q\to 0$ in the general case was
determined in [vMou]. We determine the asymptotics for $q\to 0$ in the odd case, using a different proof.

Define the mapping $\F:\cQ^{odd}\to\R^k$ by
\[
 \lb{iL01}
 \F(q)=\{\f_{2n-1}(q)\}_{1}^k,\ \ \
 \D(\l,q)=\l^{N+1}+\f_1(q)\l^{N-1}+\f_3(q)\l^{N-3}+...+\f_{N}(q).
\]
For the case $q\in \cQ^{odd}$ we denote by $d_q$ the derivative with respect to $q_1,..,q_k$. Note that $S(q)=\det d_q\F(q)$ is a polynomial, where $d_q\F$ is the
derivative of $\F$ with respect to $q$. Introduce the surface
$\cS=\{q\in\cQ^{odd}:\ \det d_q\F(q)=0\}$. Introduce  the set of
isospectral  potentials ${\rm Iso}(q)=\{p\in\cQ^{odd}:\ \
\s(q)=\s(p)\}$. Let $\# A$ denote the number of elements of the
set $A$. We formulate the main result of our paper.

\begin{theorem}
\lb{T2}
 Let $q\in\cQ^{odd}$. Then\\
i)  $\det d_q h\not=0$ iff $q\notin\cS=\{q\in\cQ^{odd}:\ \det d_q\F(q)=0\}$.\\
ii) If $q\not\in\cS$, then $h_n(q)>0$ for all $n=1,..,k$, i.e., 
all gaps are open.\\
 iii) If $\hat q_n\not=0$, for all $n=1,...,k$, then
$p\ev sq\not\in\cS$ and $\#{\rm Iso}(p)=2^k$ for any $s\in(0,\t)$
and some $\t=\t(q)>0$.
\end{theorem}

It is important to compare the two cases $q\to 0$ and $\|q\|\to
\iy$. We recall our result from [KKu2] for large $q$.

\begin{theorem}
\lb{T3} Let $q\in\cQ^{odd}$ be such that $q_i\not=q_j$ for all $i\not=j$. Then $\#{\rm Iso}(tq)=2^k k!$ for all 
$t>t_0(q)$ and some $t_0(q)>0$.
\end{theorem}

 By Theorem \ref{T2},  \ref{T3},  there exists a big difference
between the two cases: small $q$ and large $q$. Roughly speaking,
we have $\#{\rm Iso}(q)=2^k$ for small $q$ and $\#{\rm Iso}(q)=2^k
k!$ for large $q$.
Such a difference between large and small $q$ is absent in the inverse spectral theory for the continuous case.
We roughly describe the difference between the two cases: small and large $q\in\cQ^{odd}$.

If $q\to 0$, then \er{t0} gives the asymptotics
$$
\D(\l,q)=\D_M^0(\l,q) +O(\|q\|^4),\ \ as \ q\to  0, \ where \ \ \D_M^0(\l,q)=\D(\l,0)+(Af(q),\L(\l)).
$$
where  $A, f(q), \L$ are defined before Lemma 2.2 and Theorem 2.4.
Here $\D_M^0(\l,q)$ is our model Lyapunov function for small $q\in\cQ^{odd}$. If we change signs of the components $\hat q_n$ (see \er{17}),
then $\D_M(\l,q)$ is not changed. Thus if we
change signs of the components $\hat q_n$,
then the change of $\D(\l,q)$ is small. Thus we get 
$\#{\rm Iso}(q)=2^k$.

Consider $q$ large. The identity \er{23} yields the asymptotics
$$
\D(\l,tq)=\D_M(\l,tq)+O(t^{N-1}),\ \ t\to \iy,\ \ \ where \  \ \D_M(\l,q)=\prod_1^{k}(\l^2-q_n^2).
$$
Here $\D_M(\l,q)$ is our model Lyapunov function for large $q\in\cQ^{odd}$. 
Permutations and changes of signs of the components $q_n$ do not change $\D_M(\l,q)$. Thus if we change signs of the components $q_n$ and rearrange the components $q_n, n=1,..,k$,
then the change of $\D(\l,q)$ is small. Thus we get 
$\#{\rm Iso}(q)=2^kk!$ different potentials $q\in\cQ^{odd}$.
These results suggest that a global characterization of the 
isospectral set is a very hard problem.

\section {Analysis for small potentials}
\setcounter{equation}{0}

In this Sect. we prove Theorem \ref{T1}. In order to study $\D$ we
define the sets $D_j^n\ss\N^j$ of indeces by: $
D_1^n=\bigcup\limits_{j=1}^{n}{\{j\}} $ and
\[
 \lb{Sjk}
 D_j^n=\lt\{\a\in\N^j:\ 1\le\a_1<...<\a_j<n+1,
 \ \a_{s+1}-\a_s\mbox{ is odd},s=1,...,j-1\rt\},
\]
where $2\le j\le n.$ With each set $D_j^n$ we associate a
polynomial $G_j^n(\l,q), (\l,q)\in\C^{N+1}$ by
\[
\lb{G}
 G_0^n=2,\ \ G_1^n(\l,q)=\sum_{i=1}^{n}{(\l-q_i)},\ \
 G_j^n(\l,q)=\sum_{\a \in D_j^n}Q_\a(\l,q),\ 2\le j\le n\le  N+1,
\]
where $Q_\a(\l,q)=(\l-q_{\a_1})\cdot...\cdot(\l-q_{\a_j})$,
$\a=(\a_1,...,\a_j)\in\N^j$. In Sect. 4 we will prove

\begin{lemma}
\lb{2.6} 
For any $N\ge1$ the following identities are fulfilled:
\[
\lb{23}
 \D=\sum_{j=0}^{\frac {N+1}2}(-1)^{j+\frac {N+1}2} G_{2j}^{N+1},\ \
 (N \ {\rm  odd})\ \ \ {\rm and } \ \ \D=\sum_{j=0}^{\frac {N}2}
 (-1)^{j+\frac {N}2} G_{2j+1}^{N+1},\ (N \ {\rm  even}).
 \lb{Ue}\lb{Uo}
\]
\end{lemma}
\no Define two $(N+1)\ts(N+1)$ matrices $M_\n,\ \ M_\t$ by
\[
\lb{dM}
 M_\n q=(q_2,..,q_{N+1},q_1)^T,\ \ M_\t q=(q_{N+1},q_N,..,q_1)^T,\
\ q=(q_1,..,q_{N+1})^T\in\C^{N+1}.
\]
Introduce the space of homogeneous polynomials $\cP_2$ by
\[
 \cP_2=\{f:\R^{N+1}\to\R\ :\ \mbox{\ degree\ } f=2,\ f(M_\n \cdot)\ev f(\cdot)
 \},\ \ N+1=2k.
\]
Introduce $f_0$ and the vector function $f$ by
\[
 f_0(q)=(q,q),\ \ f(q)=\{f_s(q)\}_{s=1}^k,\ \ f_s(q)=(q,M_\n^s q),\ \ q\in\R^{N+1},\ \ N+1=2k.
\]
Let $\vp_m^0(\l)\ev\vp_{m}(\l,0),\ m\ge0$. Define polynomials
$\p_m^k, 1\le m\le k,$ and the $k\ts k$ matrix
$A=\{A_{n,m}\}_{n,m=1}^k$ by
\[
 \lb{9999}
 2^{-\d_{k,m}}\p_m^k(\l)=\vp_{m}^0(\l)\vp_{N+1-m}^0(\l)=
 A_{1,m}\l^{2k-2}+A_{2,m}\l^{2k-4}+...+A_{k,m}.
\]
In order to show Theorem \ref{t111} we need the following two
lemmas, proved in Sect. 4.

\begin{lemma} \lb{9990} 
The following identities are fulfilled
\[
 \lb{9989}
 \D(\l,0)=\cos(2k\arccos\frac{\l}2),\ \ \ \l\in[-2,2],\ \ \
 \l_n^0=\l_n(0)=\l_n^\pm(0)=-2\cos\frac{\pi n}{2k}, n=1,..,N,
\]
\[
 \lb{iL}
 \D(\l,0)=\sum_{m=0}^{k}
 (-1)^{m+k}R_{2m}^{2k}\l^{2m},\ \ \
 R_{2m}^{2k}=\frac{2k}{k+m}C_{2m}^{k+m},
 \ \ \ C_m^n=\frac{n!}{(n-m)!m!}.
\]
\end{lemma}
\begin{lemma}
\lb{113} i) The polynomials $f_s(q)=(q,M_\n^s q),\ s=0,...,k$ form
a basis of the space $\cP_2$.

\no ii) For any $(\l,q)\in\R\ts\cQ^{odd}$ and $1\le m\le k$ the
following identities are fulfilled
\[
 \lb{l0}
 G^{2k}_{2m}(\l,q)=R_{2m}^{2k}\l^{2m}+\sum_{n=1}^m
 g_{k+1-m,2n}(q)\l^{2m-2n},
\]
\[
 \lb{l1}
g_{m}(q)\ev g_{m,2}(q)= (-1)^{m-1}(A_{m,1}f_1(q)+...+A_{m,k}f_k(q)),
\]
for some homogeneous polynomials $g_{m,s}(q)$ with degree
$s$, and $A_{n,m}$ is given by \er{9999}.
\end{lemma}

\no Using  \er{y1} we obtain
\[
 \lb{iL00}
\D(\l,q)=\D(-\l,q)=\D(\l,-q)=\D(\l,M_\t q),\ \ \
(\l,q)\in\R\ts\cQ^{odd}.
\]
Then for any $q\in\cQ^{odd}$ the even polynomial $\D(\cdot,q)$ has
form \er{iL01}. Introduce the vectors
\[
 \lb{iL0}
 \F^0=\F(0)=(-R_{2k-2}^{2k},R_{2k-4}^{2k},...,(-1)^{k}
 R_0^{2k})^T,\ \ \ \L(\l)=(\l^{2k-2},\l^{2k-4},...,1)^T,
\]
where $\F$ is defined by \er{iL01}.

\begin{theorem}
\lb{t111} 
For any $(\l,q)\in\R\ts\cQ^{odd}$ the functions $\F$ and
$\D$ have the form
\[
 \lb{t1}
 \F(q)=\F^0+Af(q)+\F^2(q),
\]
\[
 \lb{t0}
 \D(\l,q)=\D(\l,0)+(A f(q),\L(\l))+(\F^2(q),\L(\l)),
\]
where the constant $k\ts k$ matrix $A$ and the vector function
$\F^2(q)$ satisfy
\[
 \lb{t2}
 \det A=1,\ \ \
 \|\F^2(q)\|\le K_N\|q\|^4,\ \ \ |(\F^2(q),\L(\l))|\le N
 K_N (1+|\l|^N)\|q\|^4
\]
for some  absolute constant $K_N$.
\end{theorem}
\no {\it Proof.} Substituting \er{l0}, \er{l1} into \er{Ue} and
using \er{9999} we obtain \er{t1}. Moreover, by Lemma \ref{113},
the polynomials $g_{m,2n}$ have degree $\ge4$, thus we have
$\|\F^2(q)\|\le K_N\|q\|^4$. We will show that $\det A=1$. Recall
that $\vp_n^0=\vp_n(\l,0)$ and the polynomials $\p_m^k$ are given
by
\[
 \lb{013}
 \p_m^k=\vp_{m}^0\vp_{N+1-m}^0,\ \ 1\le m\le k,\ \ N+1=2k.
\]
Then, using \er{011}, for any $1\le m<k$ we have
$$
 \p^k_{m+1}-\p^k_m=
 \vp_{m+1}^0\vp_{N-m}^0-\vp_{m}^0\vp_{N+1-m}^0=
(\l\vp_{m}^0-\vp_{m-1}^0)\vp_{N-m}^0-
\vp_{m}^0(\l\vp_{N-m}^0-\vp_{N-m-1}^0)
$$
\[
 =\vp_{m}^0\vp_{N-1-m}^0-\vp_{m-1}^0\vp_{N-1-(m-1)}^0
 =\p^{k-1}_m-\p^{k-1}_{m-1},
\]
which yields
\[
 \lb{012}
 \p^k_{m+1}-\p^k_m=\p^{k-m+1}_2-\p^{k-m+1}_1=
 \vp_2^0\vp_{2(k-m)}^0-\vp_1^0\vp_{2(k-m)+1}^0=
 \vp_{2(k-m)-1}^0.
\]
Define the matrix $B$ by
\[
 \lb{015}
 B=\{B_{n,m}\}_{n,m=1}^k,\ \ B_{n,1}=A_{n,1},\ \ B_{n,m}=
 A_{n,m}- A_{n,m-1},\ \ 2\le m\le k,\ \ 1\le n\le k.
\]
This matrix satisfies $\det B=\det A.$ We will show that $B$ is
triangular and $B_{m,m}=1,\ \ m=1,...,k$.
Using \er{9999} and \er{012} we have
$$
 \p^k_{m+1}-\p^k_m=\sum_{n=1}^{k}
 (A_{n,m+1}-A_{n,m})\l^{2(k-n)}=\sum_{n=1}^{k}B_{n,m+1}\l^{2(k-n)}
$$
$$
 =\vp_{2(k-m)-1}^0(\l)=\l^{2(k-m-1)}+O(\l^{2(k-m-2)}),\
\ 1\le m\le k-1,
$$
which yields $B_{m+1,m+1}=1,\ \ B_{n,m+1}=0,\ \ 1\le n<m+1,\ \
1\le m\le k-1$. Then $\det B=1=\det A$. Using \er{iL0}, \er{iL},
\er{t1}, we have
$$
 \D(\l,q)=\l^{N+1}+(\F(q),\L(\l))=\l^{N+1}+(\F^0,\L(\l))+(A
 f(q),\L(\l))+(\F^2(q),\L(\l)),
$$
$$
 \D(\l,0)=(\F^0,\L(\l)),\ \ \
 |(\F^2(q),\L(\l))|\le|(\F^2(q),\L(\l))|\le\|\F^2(q)\|N(1+|\l|^N),
$$
which yields \er{t0}. \BBox

Below we need 

\begin{lemma} \lb{Tiden} Let a matrix $W=\{W_{n,m}\}_{n,m=1}^k=\rt\{\cos{\pi n
m\/k}\rt\}_{n,m=1}^k$. The following identity 
%is fulfilled 
\[
 \lb{9998}
 {-2\L(\l^0_n)^T\/\D''(\l_n^0,0)}AW={2^{\d_{k,n}}\/4k}e_n
 \ \ 
\]
holds, where $e_n=\{\d_{n,m}\}_{m=1}^k,\ \ n=1,..,k.$
\end{lemma}
\no {\it Proof. } Let $\x_n={(-1)^n\/\sin^2 {\pi n\/2k}}$ and
$\p_n=\p_n^k$. For the case $q=0$ we have
\[
 \lb{9997}
 \D(\l,0)=2\cos 2kz,\ \ z=\arccos\frac\l2,\ \ \l\in[-2,2];\ \ \
 \D''(\l_n^0,0)=-2k^2\x_n,
\]
\[
 \lb{9996}
 \vp_n^0={\sin nz\/\sin z};\ \
 \p_m={\cos(2k-2m)z-\cos2kz\/2\sin^2z},\ \ \p_m(\l_n^0)=
 {\x_n\/2}\left(\cos\frac{\pi nm}k-1\right).
\]
Using \er{9996} and
$\L(\l^0_n)^TA=\{2^{-\d_{k,j}}\p_j(\l_n^0)\}_{j=1}^k$ (see
\er{9999}) we have
$$
 \left(\L(\l^0_n)^TAW\right)_{m}=
\frac12\p_k(\l_n^0)\cos\frac{\pi km}k+
\sum_{j=1}^{k-1}\p_j(\l_n^0)\cos\frac{\pi jm}k
$$
$$
 =-{\x_n\/2}\left({(\cos\pi n-1)\cos\pi
 m\/2}+\sum_{j=1}^{k-1}\cos\frac{\pi jn}{k}\cos\frac{\pi jm}{k}-
 \sum_{j=1}^{k-1}\cos\frac{\pi jm}{k}\right)=
 \d_{n,m}2^{\d_{k,n}}{k\/4}.
$$
Then, the last identity and \er{9997} yield \er{9998}. \BBox

In order to show Theorem 1.1 we need

\begin{lemma} \lb{t112} i) For any $q\in\cQ^{odd}$ the function $f=\{f_s\}_{s=1}^k,\ \ f_s(q)=(q,M_\n^s q)$
satisfies
\[
 \lb{b3}
 f(q)=W(\hat q_1^2,...,\hat q_k^2)^T,\ \ \ \ {\rm where}\ \ W=\rt\{\cos{\pi n m\/k}\rt\}_{n,m=1}^k,\ \ \det W\not=0.
\]

\no ii)  Let $\F(q),q\in\cQ^{odd}$ be given by  \er{iL01}, and let
$\hat q_n\not=0$ for all $n=1,..,k$. Then $\F:\cQ^{odd}\to\R^k$ is
a local isomorphism at any point $tq$, $t\in(0,t_0)$ for some
$t_0=t_0(q)>0$.
\end{lemma}

\no {\it Proof.} i) Recall that $N+1=2k$. Introduce vectors $\wt
e_m\in\C^{N+1}$ by
\[
 \wt e_m=\frac1{2\sqrt k}(1,s^m,...,s^{mN})^T,\ \ \
  m=1,...,N,\ \ \ s=e^{\frac{i\pi}k}.
\]
Using $s^{N+1}=1$ and \er{dM}, we have
\[
 \lb{b1}
 M_\n \wt e_m=s^m \wt e_m,\ \ \ (\wt e_l,\wt e_m)=\frac12\d_{l,m},\ \
 \ M_\t \wt e_m=s^{N+1-m} \wt e_{N+1-m},\ \ \ 1\le l,m\le
 N.
\]
Also introduce vectors $\hat e_m\in\C^{N+1}$ by
\[
 \lb{b2}
 \hat e_m=-2^{-\frac{\d_{k,m}}2}
 (s^{\frac{k+m}2}\wt e_m+s^{-\frac{k+m}2}\wt e_{N+1-m}),\
 \ \ 1\le m\le k,\ \ \ s^{\frac12}\ev e^{\frac{i\pi}{2k}}.
\]
We have $\hat e_m=(-2^{-\frac{\d_{k,m}}2})2\Re s^{\frac{k+m}2}\wt
e_m$, since $s^{N+1}=1$. Then, \er{b1} and $s^k=-1$ yield
$$
 -2^{\frac{\d_{k,m}}2}M_\t \hat e_m=s^{\frac{k+m}2}s^{N+1-m}\wt e_{N+1-m}+s^{-\frac{k+m}2}s^m\wt
 e_m=-s^{-\frac{k+m}2}\wt e_{N+1-m}-s^{\frac{k+m}2}\wt e_m=2^{\frac{\d_{k,m}}2}\hat e_m.
$$
Relations $\hat e_m\in\R^{N+1}$ and $M_\t \hat e_m=-\hat e_m$ show
that $\hat e_m\in\cQ^{odd},\ \ m=1,...,k$. Also \er{b1} yields
$(\hat e_l,\hat e_m)=\d_{l,m},\ \ 1\le l,m\le k$. This and
dim $(\cQ^{odd})=k$ show that $\hat e_m$,  $m=1,...,k$ is an
orthonormal basis in the space $\cQ^{odd}$. Define the unitary
operator $U:\cQ^{odd}\to\cQ^{odd}$ by
\[
 \lb{043}
 U(\hat e_m)=e_{m-1}-e_{N+1-m},\ \ m=1,...,k,\ \
 e_m=\{\d_{m,n}\}_{n=0}^N.
\]
The identities \er{b1}, \er{b2} yield $(\hat e_j,M_\n^n\hat e_m)=0,\ \
j\not=m,\ \ 1\le n,j,m\le k.$ Then we have
\[
 \lb{017}
 (q,M_\n^n q)=\hat q_1^2(\hat e_1,M_\n^n \hat
 e_1)+...+\hat q_k^2( \hat e_k,M_\n^n
 \hat e_k),\ \ q=\sum_{n=1}^k\hat q_n\hat e_n,
\]
and using \er{b1}, \er{b2}, we obtain
$$
 (\hat e_m,M_\n^n \hat e_m)=
 2^{-\d_{k,m}}(s^{\frac{k+m}2}\wt e_m+s^{-\frac{k+m}2}\wt e_{N+1-m}\
 ,\
 s^{\frac{k+m}2}s^{nm}\wt e_m+s^{-\frac{k+m}2}s^{n(N+1-m)}\wt
 e_{N+1-m})
$$
\[
 \lb{018}
 ={s^{-nm}+s^{nm}\/2}=\cos\frac{n m\pi}k.
\]
Thus \er{017}, \er{018} give $f(q)=W(\hat q_1^2,...,\hat
q_k^2)^T$, and \er{9998} implies $\det W\not=0$.

ii) The identities \er{t1}, \er{t2} and \er{b3} yield
$$
 \F(q)=\F^0+AW(\hat q_1^2,...,\hat q_k^2)^T+\F^2(q),\ {\rm\ where\ \ }
 \F^2(q)=O(\|q\|^4).
$$
Then for fixed $q$ we obtain
\[
 d_q\F|_{tq}=2tAW \cdot diag(\hat q_1,...,\hat q_k)
 U^T+O(t^3),\ \ t\to0,
\]
where $U$ is given by \er{043}. Due to $\hat q_n\not=0,\
n=1,...,k$ and \er{9998}, the matrix $d_q\F|_{tq}$ has an inverse
for sufficiently small $t$, since $A$, $W$ and $U$ are invertible.
\BBox

\no {\bf Proof of Theorem 1.1.} Recall that
$(-1)^{N+1-n}\D(\l_n^{\pm}(q),q)=2$ and $\D'(\l_n(q),q)=0$,
 $\l_n^0=\l_n(0)=\l_n^\pm(0)$, $n=1,...,N$ (see \er{9989}). Let
$\l_n^\pm=\l_n^\pm(q),\ \ \l_n=\l_n(q)$. Theorem \ref{t111} gives
\[
 \D(\l^0_n,0)=(-1)^{N+1-n}2=\D(\l_n^{\pm},q)=\D(\l_n^{\pm},0)+(A
f(q),\L(\l_n^{\pm}))+O(\|q\|^4),
\]
which yields
\[
 \D(\l_n^{\pm},0)-\D(\l_n^0,0)=-(A f(q),\L(\l_n^{\pm}))+O(\|q\|^4).
\]
Then, the Taylor formula implies
\[
 \lb{i1}
 \frac12\D''(\wt\l_n^{\pm},0)(\l_n^{\pm}-\l_n^0)^2=-(A
 f(q),\L(\l_n^{\pm}))+O(\|q\|^4),\ \
 \wt\l_n^{\pm}\in(\l_n^{\pm},\l_n^0).
\]
Note that by \er{t0}, $\D(\l,q)\to\D(\l,0)$ as $q\to0$, uniformly
on bounded set of $\C$ and $\l_n^{\pm}(q)\to\l_n^{\pm}(0)=\l_n^0$
as $q\to0$. Thus due to the estimate $\|f(q)\|=O(\|q\|^2)$ we
obtain
\[
 \lb{i2}
 (\l_n^{\pm}-\l_n^0)^2=\frac{-2}{\D''(\l_n^0,0)}(A
 f(q),\L(\l_n^0))+o(\|q\|^2)=O(\|q\|^2),
\]
since $\D''(\l_n^0,0)\not=0$. We will determine sharper
asymptotics. The asymptotics \er{i1} gives
$$
 {(\l_n^{\pm}-\l_n^0)^2\/2}=
\frac{-(Af(q),\L(\l_n^{\pm}))}{\D''(\wt\l_n^{\pm},0)}+O(\|q\|^4)
=\frac{-(A f(q),\L(\l^0_n))}{\D''(\l^0_n,0)}
$$
\[
 \lb{117}
-(Af(q),\L(\l^0_n))\lt(\frac1{\D''(\wt\l_n^{\pm},0)}-\frac1{\D''(\l^0_n,0)} \rt) +\frac{-(A
 f(q),\L(\l_n^{\pm})-\L(\l^0_n))}{\D''(\wt\l_n^{\pm},0)}+O(\|q\|^4),
\]
and using \er{i2}, we have
\[
\lb{115}
\left|\frac1{\D''(\wt\l_n^{\pm},0)}-\frac1{\D''(\l^0_n,0)}
\right|=O(\wt\l_n^{\pm}-\l^0_n)=O(\l_n^{\pm}-\l^0_n)=O(\|q\|),
\]
\[
 \lb{116}
 \|\L(\l_n^{\pm})-\L(\l^0_n)\|= O
 (\l_n^{\pm}-\l^0_n)=O(\|q\|).
\]
Substituting \er{115}, \er{116} into \er{117} and using \er{b3}, \er{9998} we get \er{i0}.

We will determine the asymptotics of $h_n$. Using \er{t0},
\er{t2} and \er{b3}, \er{9998} we get
$$
 \D(\l_n^0,q)=\D(\l_n^0,0)+(A f(q),\L(\l_n^0))+O(\|q\|^4)=
 (-1)^n2+(AW(\hat q_1^2,...,\hat q_k^2)^T,\L(\l_n^0))+O(\|q\|^4)
$$
\[
 \lb{9979}
 =(-1)^n2+\L(\l_n^0)^TAW(\hat q_1^2,...,\hat q_k^2)^T+O(\|q\|^4)=
 (-1)^n\left(2+\frac{2^{\d_{k,n}}k}{4\sin^2\frac{\pi n}{2k}}
 \hat q_n^2\right)+O(\|q\|^4),
\]
and, for small $t$, the Taylor formula implies
\[
 \lb{9978} \D(\l^0_n+t,q)=\D(\l^0_n,q)+{t^2\/2}\D''(\wt\l^0_n(t),0)+tO(\|q\|^2)+O(\|q\|^4),\ \ |\wt\l^0_n(t)-\l^0_n|\le |t|.
\]
Let $t_0=C\|q\|^2$. Then, for sufficiently large $C>0$ and
sufficiently small $\|q\|$, the identities
$\sign\D(\l^0_n,0)=-\sign\D''(\l^0_n,0)$ and \er{9978} yield
\[
 |\D(\l^0_n\pm t_0,q)|<|\D(\l^0_n,q)|,\ \ \
 \D(\l^0_n+\t,q)=\D(\l_n^0,q)+O(\|q\|^4),\ \ \ |\t-\l^0_n|\le t_0,
\]
which gives
\[
 \lb{9977}
 \D(\l_n,q)=\D(\l_n^0,q)+O(\|q\|^4)=(-1)^n\left(2+
 \frac{2^{\d_{k,n}}k}{4\sin^2\frac{\pi n}{2k}}
 \hat q_n^2\right)+O(\|q\|^4).
\]
Then, using \er{9977}, \er{1d} and ${\rm
arccosh}^2(1+x)=2x+O(x^2)$, we get \er{9976}. \BBox

\section {Proof of Theorem \ref{T2}}
\setcounter{equation}{0}

In this section $q\in\cQ^{odd}$. Introduce the matrix
\[
\lb{4V}
V(a)=\left(\begin{array}{cccc}
{a_1^{2(k-1)}} & {a_1^{2(k-2)}} & ... & {1} \\
{a_2^{2(k-1)}} & {a_2^{2(k-2)}} & ... & {1} \\
... & ... & ... & ... \\
{a_N^{2(k-1)}} & {a_N^{2(k-2)}} & ... & {1}
\end{array}\right),\ \ \ a=(a_1,...,a_k)^T\in \R^k.
\]
Define the mapping $ H:\cQ\to \R^k$ by $q\to  H(q)={\{
H_n(q)\}}_{1}^k,$ with the components
\[
 \lb{061}
  H_n= \D(\l_n,q)=2(-1)^{N+1-n}\cosh h_n.
\]
\begin{lemma} \lb{T41} i) The mappings $ H:\cQ^{odd}\to\R^k$ and
$\wt\l=\{\l_n\}_{n=1}^k:\cQ^{odd}\to\R^k$ are real analytic and
\[
\lb{4p}
 d_q H(q)=V(\wt\l(q))d_q \F(q).
\]
\no ii) The mapping $ h:\cQ^{odd}\to\R^k$ is continuous and the
mapping $ h:\wt\cQ^{odd}\to\R^k$ is real analytic, where
$\wt\cQ^{odd}=\{q\in\cQ: h_m(q)>0,\ m=1,...,k\}$.

\end{lemma}

 \no {\it Proof.} i) 
 The function $\D'(\l,q)$ is a polynomial of
degree $N$ in $\l$, whose coefficients are polynomials of $q$. Therefore, its roots $\l_n, n=1,..., N$ are
continuous functions of $q$. These roots are simple,
thus they are real analytic on $\cQ^{odd}$.
Then the function
$H_n(q)=\D(\l_n(q),q), q\in \cQ^{odd}$ is real analytic,
since $\D(\l,q)$ is a polynomial. 
The identities ${\D'}(\l_n(q),q)=0$ and

$$
 \pa_m H_n=\pa_m \lt(\D(\l_n(q),q)\rt)=
\l_n^{N-1}\pa_m \f_1+\l_n^{N-2}\pa_m \f_2+...+\pa_m
\f_N+\D'(\l_n(q),q)\pa_m\wt\l_n(q)
$$
yield \er{4p}.

\no ii) Note, that arccosh is a continuous function on the set
$\{x\in\R,\ x\ge1\}$  and it is real analytic on the set $\{x\in\R,\
x>1\}$. Then, $ h_n(q)=$ arccosh $((-1)^{N+1-n} H_n(q)/2)$ is a
continuous mapping on $\cQ^{odd}$ and real analytic mapping on the
open set $\wt \cQ^{odd}$, since
 we have $(-1)^{N+1-n} H_n(q)\ge2, q\in\cQ^{odd}$ and $(-1)^{N+1-n} H_n(q)>2,\ q\in\wt\cQ^{odd}$.
\BBox.

Introduce the set of indexes
\[
 \cN=\{\n=\{\n_n\}_{1}^k:\ \ \n_n\in\{-1,+1\}\ \},\ \ \#\cN=2^k.
\]
For $(\n,\ve)\in\cN\ts\R_+$ define the sets
\[
 \cZ^\n_\ve=\{q\in\cQ^{odd}:\ \ \|q\|<1,\ \ \n_n\hat q_n>\ve,\ n=1,...,k\},
\]
\[
 \cZ^\n(\ve,t)=\{q:\ \ q=\t p,\ \ p\in\cZ^\n_\ve,\ \
 \t\in(0,t)\}.
\]
Note that $\cZ^\n_\ve\not=\es$ for
$0<\ve<\frac1{\sqrt{N+1}}$.

\begin{lemma} \lb{035} 
For any $t\in(0,t_0(\ve))$ and $0<\ve<\frac1{\sqrt{N+1}}$,
the mapping $\F:\ \cZ^\n(\ve,t)\to\R^k$ is an injection and a
local isomorphism for some $t_0(\ve)$, which depends only on
$\ve$. \lb{038}
\end{lemma}

\no {\it Proof.} Let $\wt p,\wt q\in\cZ^\n(\ve,t)$ and let $\wt
p\not=\wt q$. Then $\wt q=tq,\ \ q\in\cZ^\n_\ve$ and $\wt
p=tq+t\d,\ \ q+\d\in\cZ^\n_\ve$ and $0<\|\d\|<1$. Using \er{t1}
and a polynomial $\F^2(q)=O(\|q\|^4)$, we obtain
$$
 \F(tq+t\d)-\F(tq)=tAWdiag(2\hat
q_n+\hat\d_n)(t\hat\d_1,...,t\hat\d_k)^T+r(t)(t\hat\d_1,...,t\hat\d_k)^T
$$
\[
 =t^2(AWdiag(2\hat
 q_n+\hat\d_n)+O(t^2))(\hat\d_1,...,\hat\d_k)^T\not=0,
\]
since $|2\hat q_n+\hat\d_n|>2\ve$ and $AWdiag(2\hat
 q_n+\hat\d_n)$ is invertible . Then for sufficiently
small $t$ we have $\F(\wt p)\not=\F(\wt q)$. 
 By Theorem 1.2, $\F$ is a local isomorphism. \BBox

{\no \bf Proof of Theorem \ref{T2}.} Suppose $h_n(q)=0$, for some
$q\in\cQ^{odd}$ and some $n$. Then $| H_n(q)|=2$. If
$q\not\in\cS=\{q\in\cQ^{odd}:\ \det d_q\F(q)=0\}$, then $ H$ is a
local isomorphism at the point $q$ (see \er{4p}). Thus there
exists some point $p\in\cQ^{odd}$ such that $| H_n(p)|<2$, which contradicts $| H_n(p)|\ge2,\ \ p\in\cQ^{odd}$. We obtain $q\in\cS$.

Using \er{061}, \er{4p}, we have that $ h$ is a local isomorphism
at any point $q\not\in\cS$, since $V(\wt\l(q))$ is the Vandermond
matrix and $\l_n(q)\not=\l_m(q),\ \ \ m\not=n.$

Let $q\in\cQ^{odd},\ \ \|q\|=1$ and let $|\hat q_n|>\ve$ for all
$n=1,...,k$. Let $0<2t<t_0({\ve\/3})$, where the function
$t_0(\cdot)$ is defined in Lemma \ref{035}. Let $\wt q=tq$. Then
the first statement of our theorem yields $\wt q\not\in\cS$. 

Consider the first case $\wt
p\not\in\bigcup\limits_\n\cZ^\n({\ve\/3},2t)$ and $\|\wt p\|<2t$.
Then $\wt p=2tp,\ \ \|p\|<1,\ \ |p_n|<{\ve\/3}$ for some $n$,
which gives $|\hat q_n^2-4\hat p_n^2|\ge\frac{5\ve^2}9$. Thus for
sufficiently small $t$ the identity \er{t1} yields
$$
 \|\F(\wt q)-\F(\wt p)\|=\|\F(tq)-\F(2tp)\|=t^2AW(\hat q_1^2-4\hat p_1^2,...,\hat q_k^2-4\hat
 p_k^2)^T+O(t^4)
$$
\[
 \lb{031}
 \ge C t^2 |\hat q_n^2-4\hat p_n^2|\ge
 C t^2 \frac{5\ve^2}9,\ \ C={\|W^{-1}A^{-1}\|\/2}.
\]

Consider the second case $\wt
p\not\in\bigcup\limits_\n\cZ^\n({\ve\/3},2t)$ and $\wt p=2tp,\ \
\|p\|\ge1$ for some $p\in\cQ^{odd}$. Then, using
$\f_1(q)=\frac{\|q\|^2}2+N+1$, we obtain
\[
 \lb{032}
 \|\F(\wt q)-\F(\wt p)\|=\|\F(tq)-\F(2tp)\|\ge\frac32t^2\ge t^2.
\]
The estimates \er{031}, \er{032} and the fact that $\F$ is a local
isomorphism in the sets $\cZ^\n({\ve\/3},2t),\ \ \n\in\cN$ (see
Lemma \ref{038}) yield
\[
 \lb{033}
 dist(\F(\wt q),\pa\F(\cZ^\n({\ve\/3},2t)))=dist(\F(\wt q),\F(\pa\cZ^\n({\ve\/3},2t)))\ge t^2
 \min\left(C\frac{5\ve^2}9,1\right),
\]
for all $\n\in\cN$.

For any $\n\in\cN$ we take $2^k$ points $\wt q^\n=tq^\n=t(\n_1\hat
q_1,...,\n_k \hat q_k)\in\cZ^\n({\ve\/3},2t)$. Note that $\wt
q=\wt q^{\n_0},\ \ \n_0=(1,...,1)\in\cN$. Then \er{t1} yields
\[
 \lb{034}
 dist(\F(\wt q^\n),\F(\wt q))=\|\F(\wt q^\n)-\F(\wt
 q)\|=\|\F(tq^\n)-\F(tq)\|=O(t^4),\ \ {\rm for\ all\ \n\in\cN}.
\]
Using \er{033}, \er{034} and the fact that $\F$ is a local
isomorphism in the sets $\cZ^\n({\ve\/3},2t),\ \ \n\in\cN$, we
have
\[
 \F(\wt q)\in \F(\cZ^\n({\ve\/3},2t)),{\rm\ \ for\ all}\ \ \n\in\cN,
\]
since $dist(\F(\wt q),\F(\wt q^\n))<dist(\F(\wt
q),\pa\F(\cZ^\n({\ve\/3},2t)))$ and $\F(\wt
q^\n)\in\F(\cZ^\n({\ve\/3},2t))$ for all $\n\in\cN$ and for
sufficiently small $t$. Then there exist $2^k$ distinct points
$v^\n\in\cZ^\n({\ve\/3},2t)$ such that $\F(\wt q)=\F(v^\n),\ \
\n\in\cN$.

We show that there are no other points, i.e. we have exactly $2^k$
distinct points $v^\n$. Suppose $\F(\wt q)=\F(u)$ for some
$u\in\cQ^{odd}$, $u\not=v^\n$ for all $\n\in\cN$. Then
$u\not\in\bigcup_\n\cZ^\n({\ve\/3},2t)$,  since by Lemma \ref{035}, $\F$ is an injection in the sets $\cZ^\ve({\ve\/3},2t)$.
Using $\f_1(q)=\frac{\|q\|^2}2+N+1$, we have $\|u\|=\|\wt q\|=t<2t$.
Thus \er{031} holds for $\wt q, u$, which contradicts $\F(\wt q)=\F(u)$. \BBox

\section {Appendix}
\setcounter{equation}{0}

In this Section we prove Lemmas \ref{2.6}-\ref{113}, i.e., we determine some identities for the Lyapunov
function $\D(\l,q)$ for the general case $q\in\cQ$, not only
$q\in\cQ^{odd}$. The results of this Section are used also to study the case of $q$ large  in [KKu2]. Introduce the sets $T_j^n\ss\N^j$ by
\[
 T_j^n=\rt\{\a=\{\a_s\}_{1}^j:\ 1\le \a_1<...<\a_j<n+1\ev \a_{j+1},
 \a_{s+1}-\a_s\mbox{ is odd, } s=1,...,j \rt\},
\lb{2a}
\]
$1 \le j \le n$.  We need the following simple properties of the sets $T_j^n$.

\begin{lemma}
\lb{2.1}  Let $E_{j}^{n+1}=\{\a=(\r,n+2):\r\in T_{j-1}^{n+1}\}$.
The following relations are fulfilled
\[
 T_j^n\ss T_j^{n+2},\ \ \ {\rm any}\ \ 1\le j\le n,
\lb{2b}
\]
\[
\lb{2c} T_{j}^{n+2}=T_{j}^n\cup E_{j}^{n+1} ,\ \ \ T_{j}^n\bigcap
E_{j}^{n+1}=\es ,\ \ \ \ {\rm any} \ \ 2\le j\le n.
\]
\end{lemma}

{\no \it Proof.} In order to show \er{2b} we note that if $
n+1-\a_j$ is odd, then $ n+3-\a_j$ is odd also. The definition
\er{2a} implies \er{2b}. We prove \er{2c}. The definition \er{2a}
provides $ \{\a=(\r,n+2):\r\in T_{j-1}^{n+1}\} \ss T_{j}^{n+2} $.
\er{2b} gives $T_{j}^n \ss T_{j}^{n+2} $, then we have $
T_{j}^n\cup E_{j}^{n+1}\ss T_{j}^{n+2}. $ We show the opposite
inclusion. Suppose that
\[
\a=(\a_1,...,\a_{j-1}),\ \ \  (\a,\a_{j})\in T_{j}^{n+2}. \lb{2d}
\]
If $ \a_{j}=n+2$, then $ n+2-\a_{j-1} $ is odd and we have $
 \a\in T_{j}^{n+1}
$.
 Let $\a_{j}<n+2$. If $n+1-\a_{j}$ is odd, then
$(\a,\a_{j})\in T_{j}^n$. If $n+1-\a_{j}$ is even, then
$n+3-\a_{j}$ is even and \er{2d} is not fulfilled. The opposite
inclusion is proved and we have the identity \er{2c}.

\no The sets $T_{j}^n$ and $E_{j}^{n+1}$ are disjoint, since by
definition \er{2a}, the element $(\a,n+1)\not\in T_{j}^n$. \BBox

Recall that for $q\in\cQ$ and the multi-index $\a$  the polynomial
$Q_\a(\l,q)$ is given by
\[
\lb{2e} Q_\a(\l,q)=(\l-q_{\a_1})\cdot...\cdot(\l-q_{\a_j}),\ \ \ \
\a=(\a_1,..,\a_j)\in\N^j .
\]
With each set $T_j^n$ and $q\in\cQ$ we associate the polynomial
$F_j^n(\l,q)$  by
\[
 F_0^n\ev 1,\ \ \ \
 F_j^n(\l,q)=\sum\limits_{\a\in T_j^n} Q_\a(\l,q)
 ,\ \ \ \ \ \ \ {\rm if}\ \ \ \ \  1\le j\le n\le N+1.
 \lb{2f}
\]
In order to prove Lemma \ref{2.3}
we need some properties of $F_j^n$.

\begin{lemma}
\lb{2.2} For any $(\l,q)\in\C^{N+1}$ the following identities are
fulfilled :
\[
 F_{j+1}^{n+2}(\l,q)=(\l-q_{n+2})F_j^{n+1}(\l,q)+F_{j+1}^n(\l,q),\ \ \ 1\le j \le n-1,
 \lb{2g}
\]
\[
 F_{n+1}^{n+1}(\l,q)=(\l-q_{n+1})F_n^n(\l,q),\ \ \ n\ge1,
 \lb{2h}
\]
\[
 F_n^{n+1}(\l,q)=(\l-q_{n+1})F_{n-1}^n(\l,q),\ \ \ n\ge2.
 \lb{2i}
\]
\end{lemma}

{\no \it Proof.} Let $F_j^n=F_j^n(\l,q)$ and $Q_\a=Q_\a(\l,q)$.
Using definition \er{2f} and \er{2b} we obtain
$$
 (\l-q_{n+2})F_j^{n+1}+F_{j+1}^n
 =\sum_{\a\in T_j^{n+1}}(\l-q_{n+2})Q_{\a}
 +\sum\limits_{(\a,\a_{j+1})\in T_{j+1}^n}
 (\l-q_{\a_{j+1}})Q_{\a}
 $$
 $$
 =\sum_{(\a,\a_{j+1})\in E_{j+1}^{n+1} \cup T_{j+1}^n}
(\l-q_{\a_{j+1}})Q_{\a} =\sum\limits_{\a\in
T_{j+1}^{n+2}}Q_{\a}=F_{j+1}^{n+2},
$$
which implies \er{2g}. Definition \er{2f} gives 
$ F_{n+1}^{n+1}=(\l-q_{n+1})\cdot...\cdot(\l-q_1)=
 (\l-q_{n+1})F_n^n$ and 
 $ F_n^{n+1}=(\l-q_{n+1})\cdot...\cdot(\l-q_2)=
 (\l-q_{n+1})F_{n-1}^n,
$ which yield \er{2h} and \er{2i}. \BBox

We prove some identities for the polynomials $\vp_n$ and $\vt_n$.

\begin{lemma}
\lb{2.3} For any $n\ge 1$ the following identities are fulfilled:
\[
\vp_{n+1}=(-1)^{n\/2}\sum_{j=0}^{n\/2}(-1)^j F_{2j}^n,\ \ \ \ \
 \vt_{n+1}=(-1)^{n\/2}\sum_{j=0}^{\frac {n-2}2}(-1)^j F_{2j+1}^n,\ \ \  \ \ n\ {\rm even,}
 \lb{Pev}\lb{Qev}
\]
\[
\vp_{n+1}=(-1)^{n-1\/2}\sum_{j=0}^{n-1\/2}(-1)^jF_{2j+1}^n,\ \ \ \
\vt_{n+1}=(-1)^{\frac {n+1}2}\sum_{j=0}^{\frac {n-1}2}(-1)^j
F_{2j}^n,\ \ \ \ \ n\ {\rm odd}.
 \lb{Pod}\lb{Qod}
\]
\end{lemma}

{\no \it Proof.} We show \er{Pev} and \er{Qod} for $\vp_n$ by
induction. The proof for $\vt_n$ is similar. Using \er{1b} we have
$$
 \vp_{0}=0 ;\ \vp_1=1; \ \  \ \vp_2(\l,q)=(\l-q_1);
 \ \ \vp_3(\l,q)=(\l-q_2)(\l-q_1)-1;
$$
$$
 \vp_4(\l,q)=(\l-q_3)(\l-q_2)(\l-q_1)-(\l-q_3)-(\l-q_1).
$$
By induction, suppose that identities \er{Pev} and \er{Pod}
are fulfilled from $1$ to $n\ge2$. Let $\vp_n=\vp_{n}(\l,q),
F_j^{n}=F_j^{n}(\l,q)$. If $n$ is odd, then from Lemma \ref{2.2}
we obtain
$$
\vp_{n+2}=(\l-q_{n+1})\vp_{n+1}-\vp_n=(\l-q_{n+1})(-1)^{\frac
{n-1}2}\sum_{j=0}^{\frac {n-1}2}(-1)^j F_{2j+1}^{n}
 -(-1)^{\frac {n-1}2}\sum_{j=0}^{\frac {n-1}2}(-1)^jF_{2j}^{n-1}
$$
$$
 =(\l-q_{n+1})F_{n}^{n}+(-1)^{\frac {n-1}2}
 \sum_{j=0}^{\frac {n-1}2-1}{(-1)^j(\l-q_{n+1})F_{2j+1}^{n}}
 -(-1)^{\frac {n-1}2}-(-1)^{\frac {n-1}2}
 \sum_{j=1}^{\frac {n-1}2}{(-1)^j F_{2j}^{n-1}}
 $$
$$
= F_{n+1}^{n+1}+(-1)^{\frac {n+1}2}F_0^{n+1}
 -(-1)^{\frac {n-1}2}\sum_{j=1}^{\frac {n-1}2}
 {(-1)^j((\l-q_{n+1})F_{2j-1}^{n}+F_{2j}^{n-1})}
 $$
$$
=(-1)^{\frac{n+1}2}((-1)^{\frac {n+1}2}F_{n+1}^{n+1}
 +(-1)^0F_0^{n+1}+\sum_{j=1}^{\frac {n-1}2}{(-1)^jF_{2j}^{n+1}})=
 (-1)^{\frac {n+1}2}\sum_{j=0}^{\frac {n+1}2}{(-1)^j
 F_{2j}^{n+1}}.
$$
\no If $n$ is even, then  Lemma \ref{2.2} yields
$$
 \vp_{n+2}=(\l-q_{n+1})\vp_{n+1}-\vp_n=(\l-q_{n+1})(-1)^{\frac {n}2}
 \sum_{j=0}^{\frac {n}2}{(-1)^j F_{2j}^{n}}
-(-1)^{\frac {n-2}2}\sum_{j=0}^{\frac {n-2}2}{(-1)^j
 F_{2j+1}^{n-1}}
$$
$$=
 (\l-q_{n+1})F_{n}^{n}+(-1)^{\frac {n}2}\sum_{j=0}^
 {\frac {n-2}2}{(-1)^j((\l-q_{n})F_{2j}^{n}+F_{2j+1}^{n-1})}
$$
$$
= F_{n+1}^{n+1} +(-1)^{\frac {n}2}\sum_{j=0}^
 {\frac {n-2}2}{(-1)^j F_{2j+1}^{n+1}}=
 (-1)^{\frac {n}2}\sum_{j=0}^{\frac {n}2}{(-1)^j F_{2j+1}^{n+1}},
$$
which gives the identities \er{Pev} and \er{Pod} for any
$n\ge3$.\BBox

We need some properties of the sets $D_j^n$ and the polynomials
$G_j^n$ given by \er{Sjk} and \er{G}.

\begin{lemma}
\lb{2.5}   For any $1\le j<n$ the following identities are
fulfilled:
\[
\lb{b311}
 D_n^n=T_n^n,\ \mbox{for $n\ge1$};\ \ \
 D_j^n=T_j^n\cup T_j^{n-1},\ T_j^n\cap T_j^{n-1}=\es,
\]
\[
\lb{b4}
 G_n^n(\l,q)=F_n^n(\l,q),\ \mbox{for $n\ge1$};\ \ \
 G_j^n(\l,q)=F_j^n(\l,q)+F_j^{n-1}(\l,q).
\]
\end{lemma}

{\no \it  Proof.} The definitions \er{Sjk}, \er{2a} of $D_j^n,\
T_j^n$ give the first identities in \er{b311}, \er{b4} and
$T_j^n\cup T_j^{n-1}\ss D_j^n$. We will prove the opposite
inclusion. Suppose that $\a\in D_j^n$. We have two cases. Firstly,
if $n+1-\a_j$ is odd, then by definition \er{2a}, $\a\in T_j^n$.
Secondly, if $n+1-\a_j$ is even, then we have $n>\a_j$ and
$n-\a_j$ is odd. The inclusion $\a\in T_j^{n-1}$ is fulfilled.
Also we have showed that sets $T_j^n$ and $T_j^{n-1}$ are
disjoint.\BBox

{\no \bf Proof of Lemma \ref{2.6}.} Let $\D=\D(\l,q),
\vp_n=\vp_n(\l,q)....$ If $N$ is odd, then using Lemma 2.3-2.5, we
have
$$
 \D=\vp_{N+2}+\vt_{N+1}=(-1)^{\frac {N+1}2}
 \sum_{j=0}^{\frac {N+1}2}{(-1)^j F_{2j}^{N+1}}
 +(-1)^{\frac {N+1}2}\sum_{j=0}^
 {\frac {N-1}2}{(-1)^j F_{2j}^{N}}
$$
$$
 =(-1)^{\frac {N+1}2}(-1)^{\frac {N+1}2}F_{N+1}^{N+1}
 +(-1)^{\frac {N+1}2}\sum_{j=0}^{\frac {N-1}2}
 {(-1)^j(F_{2j}^{N+1}+F_{2j}^{N})}=
 (-1)^{\frac {N+1}2}\sum_{j=0}^{\frac {N+1}2}{G_{2j}^{N+1}}.
$$
If $N$ is even, then similar arguments yield
$$
 \D=\vp_{N+2}+\vt_{N+1}=(-1)^
 {\frac {N}2}\sum_{j=0}^{\frac {N}2}
 {(-1)^j F_{2j+1}^{N+1}(\l,q)}
 +(-1)^{\frac {N}2}\sum_{j=0}^{\frac {N-2}2}
 {(-1)^j F_{2j+1}^{N}(\l,q)}
$$
$$
=(-1)^{\frac {N}2}
 (-1)^{\frac {N}2}F_{N+1}^{N+1}
 +(-1)^{\frac {N}2}\sum_{j=0}^{\frac {N-2}2}
 {(-1)^j(F_{2j+1}^{N+1}+F_{2j+1}^{N})}=
 (-1)^{\frac {N}2}\sum_{j=0}^{\frac {N}2}
 {(-1)^j G_{2j+1}^{N+1}}.\BBox
$$

\begin{lemma}
For any $(\l,q)\in\C^{N+1}$
 the following identities are fulfilled
\[
 \lb{y1}
 G^{N+1}_j(\l,M_{\o}q)=G^{N+1}_j(\l,q),\ \ N-j\ \ is\ \ odd;
 \ \ \D(\l,M_{\o}q)\equiv  \D(\l,q),\ \ \ \o\in\{\n,\t\}.
\]
\end{lemma}

\no {\it Proof.} Define the functions $\n,\t,
e:\{1,...,N+1\}\to\{1,...,N+1\}$ by
\[
 \lb{Sjk00}
 \n(N+1)=1,\ \ \n(i)=i+1,\ \ 1\le i\le N,\ \ \t(i)=N+2-i,\ \ e(i)=i,\ \ 1\le i\le  N+1.
\]
In order to prove $G_j^{N+1}(\l,M_{\o}q)\ev G_j^{N+1}(\l,q)$ we
need to show $\o D_j^{N+1}=D_j^{N+1}$, where $\o
D_j^{N+1}=\{\b\in\N^j:\ \ \b_i=\o(\a_i),\ \ \a\in D_j^{N+1},\ \
1\le i\le j\}$, since we have \er{G}. Suppose $j=2m$ and $N+1=2k$
are even. The proof in the odd case is similar. Definitions
\er{Sjk}, \er{Sjk00}  yield $\o D_{2m}^{2k}\ss
D_{2m}^{2k},\ \ \o\in\{\n,\t\}.$ Also, using $\o^2=e,\ \
\o\in\{\n,\t\}$, we have $D_{2m}^{2k}=\o^2 D_{2m}^{2k}\ss \o
D_{2m}^{2k}$ and then $\o D_{2m}^{2k}=D_{2m}^{2k}$.  Using
\er{Ue}, we obtain $\D(\l,M_{\o}q)\equiv  \D(\l,q)$. \BBox

{\no \bf Proof of Lemma \ref{9990}.}  For the functions
$\vp^0_m(\l)=\vp_m(\l,0)$ and $\vt^0_m(\l)=\vt_m(\l,0)$ equation
\er{1b} yields
\[
 \lb{011}
 \vp^0_{0}=0,\ \vp^0_1=1,\ \vp^0_{n+1}=\l\vp^0_n-\vp^0_{n-1},\
 \vt^0_{0}=1,\
 \vt^0_1=0,\ \vt^0_{n+1}=\l\vt^0_n-\vt^0_{n-1}.
\]
Using $\D^0_{n+1}\ev\vp^0_{n+2}+\vt^0_{n+1}$, we have
$$
 \D^0_0=2,\ \D^0_1=\l,\ \D^0_2=\l^2-2,\ \D^0_{n+1}=\l\D^0_n-\D^0_{n-1}.
$$
Using similar recurrence formulas for the Chebyshev polynomials
$T_n(\l)=\cos(n\arccos(\l))$
\[
 \lb{009}
 T_0(\l)=1,\ T_1(\l)=\l,\ T_{n+1}(\l)=2\l T_n(\l)-T_{n-1}(\l)
\]
(see [AS]) we obtain $\D_{n+1}(\l,0)=2T_{n+1}({\l\/2})$. Thus we
have \er{9989} and \er{iL}, see [AS]. \BBox

{\no \bf Proof of Lemma \ref{113}.} i) Firstly, $f_s\in\cP_2$,
since $(M_\n q,M_\n q)=(q,q)$ and $f_s(M_\n q)=(M_\n
q,M_\n^{s+1}q)=(q,M_\n^s q)=f_s(q)$. Secondly, let $f\in\cP_2$.
Then $ f(q)=\sum\limits_{1\le l,m\le 2k} c(l,m)q_l q_m,$ for some
$c(l,m)=c(m,l)\in\R $. Define $q_{l+2kj}=q_{l},\ \
c(l+2kj,m+2ks)=c(l,m),\ \ 1\le l,m\le 2k,\ \ s,j\in\Z$. Then using
the identity $M_{\n}^j q=(q_{j+1},...,q_{2k+j}),\ \ \ j\in\Z$ we
obtain
\[
 f(M_{\n}^j q)=\sum\limits_{1\le l,m\le 2k} c(l,m) q_{l+j} q_{m+j},\ \ j\in\Z.
\]
This and $f(\cdot)\ev f(M_{\n}^j\cdot)$ yield $c(l,m)=c(l+j,m+j)$.
Then we have
\[
 c(1,m)=c(j+1,m+j),\ \ 1\le m,j\le 2k,
\]
which gives
\[
 f(q)=\sum_{m=1}^{2k}c(1,m)\sum_{j=1}^{2k}q_j q_{m+j}=\sum_{m=1}^{2k}c(1,m)(q,M_{\n}^m q).
 \lb{i311}
\]
Furthermore,
$(q,M_\n^{k+j}q)=(M_\n^{k-j}q,M_\n^{2k}q)=(M_\n^{k-j}q,q)=f_{k-j}(q)$.
Using this and \er{i311} and linear independence of $f_s$ we see
that $f_s,\ \ s=0,...,k$ is a basis in the space $\cP_2$.

ii) The identities \er{G},\er{Ue},\er{iL} and
$G_{2m}^{2k}(-\l,\cdot)=G_{2m}^{2k}(\l,\cdot),\ \l\in\R$ give
\er{l0}. Moreover  $G^{2k}_{2m}(\cdot,M_\n q)\ev
G^{2k}_{2m}(\cdot,q),\ q\in\cQ^{odd}$ give $g_m(M_\n q)=g_m(q),\
q\in\cQ^{odd}$, which yields $g_{m}\in \cP_2$. The proof of \er{l1}
needs Lemma \ref{044} and will be given below. \BBox

\no In order to show \er{l1} we prove Lemma \ref{044}. For $0\le
r<n-m$ introduce the set
\[
 \lb{111}
 S^{m,n}_{r}=\{\a=\{\a_s\}_{s=0}^{r+1}\in\Z^{r+2}:\ m=\a_0<...<\a_{r+1}=n,\ \a_{s+1}-\a_s\ is\ odd,\ 0\le s\le
 r\}.
\]
We see that $\# S^{m,n}_r$ depends only on $r$ and $n-m$, that is
$\#S^{m+l,n+l}_r=\#S^{m,n}_r,\ l\in\Z$. Also we see that if
$n-m-1-r$ is odd, then $S^{m,n}_{r}=\varnothing$, since we have
the identity
$$
 n-m=\sum_{i=0}^r\b_i,\ \ \ \b_i=\a_{i+1}-\a_i\ is\ odd.
$$
In order to prove \er{l1} we need

\begin{lemma}
\lb{044} 
Let $n-m-r$ be odd for some $0\le r<n-m$. Then the
following identities hold
\[
 \lb{s2}
 \#S^{m,n}_{0}=1;\ \ \ \ \ \#S^{m,n}_{r}=\#T_r^{n-m-1},\
 r\ge1.
\]
\end{lemma}
\no {\it Proof.} The set $S^{m,n}_0$ consists of one element
$(m,n)\in\Z^2$. Suppose $r\ge1$. Using the identity
$\#S^{m+l,n+l}_r=\#S^{m,n}_r,\ l\in\Z$, we need to show
$\#S^{0,n-m}_{r}=\#T_r^{n-m-1}$. Introduce the set
\[
 \wt S^{0,n-m}_{r}=\{(0,\a,n-m)\in\Z^{r+2}:\ \a\in
 T_r^{n-m-1}\}.
\]
The definition of $T$ and $S$ (see \er{2a}, \er{111}) give
$S^{0,n-m}_{r}\ss\wt S^{0,n-m}_{r}$. In order to prove $\wt
S^{0,n-m}_{r}\ss S^{0,n-m}_{r}$ we need to show that $\a_1-0$ is
odd for $\a\in T_r^{n-m-1}$. From the condition $n-m-r$ is odd and
the identity
\[
 n-m=\a_1-0+\sum_{i=1}^r\b_i,\ \ \ \b_i=\a_{i+1}-\a_i\ is\ odd\
\]
we obtain $\a_1-0$ is odd for $\a\in T_r^{n-m-1}$. Thus we have
$S^{0,n-m}_{r}=\wt S^{0,n-m}_{r}$ and then $\#S^{0,n-m}_{r}=\#\wt
S^{0,n-m}_{r}=\# T_r^{n-m-1}$. 
\BBox

{\no \bf Proof of \er{l1}.}
 We fix $m$. Introduce the sets ( $D_{j}^{k}$ is
defined by \er{Sjk} )
\[
 \lb{p1}
 P^{2m}_{j}=\{\a=\{\a_s\}_1^{2m}\in D^{2k}_{2m}:\ \a_1=1,\ \a_s=j+1\ for\ some\ s\},\ \ \ 1\le j\le k.
\]
Recall that $\d_{k,j}$ is the Kroneker symbol.
Using $g_m\in\cP_2$ and Lemma \ref{113},i) and \er{G}, \er{l0} we
obtain
\[
 g_{m}(q)=a_{m,1}f_1(q)+...+a_{m,k}f_k(q),\ \ \
 a_{k+1-m,j}=2^{-\d_{k,j}}\#P^{2m}_j,\ 1\le j\le n.
\]
 The following decompositions
are fulfilled
\[
 \lb{p2}
 P^{2m}_{j}=\bigcup_{i=2}^{j+1}P^{2m}_{i,j},\ \ \ j<2m;\ \ \ \ \
 P^{2m}_{j}=\bigcup_{i=2}^{2m}P^{2m}_{i,j},\ \ \ 2m\le j,
\]
where $ P^{2m}_{i,j}=\{\a\in D^{2k}_{2m}:\ \a_1=1,\ \a_i=j+1 \}$.
We see that the set $P^{2m}_{j}$ is a union of disjoint sets. The
following identity
\[
 \#P^{2m}_{i,j}=(\#S^{1,j+1}_{i-2})\cdot(\#S^{j+1,2k+1}_{2m-i})
\]
and  \er{s2} gives
\[
 \lb{112}
 \#P^{2m}_{i,j}=0,\ \ \ j-i\ is\ even,\ \ \
 \#P^{2m}_{i,j}=\#T^{j-1}_{i-2}\cdot\#T^{2n-j-1}_{2m-i},\ \ \ j-i\ is\  odd.
\]
Thus \er{p2}, \er{112} yield
$$ \#P^{2m}_{j}=\#T_{j-1}^{j-1}\cdot\#T_{2m-2-(j-1)}^{2k-2-(j-1)}+\#T_{j-3}^{j-1}\cdot\#T_{2m-2-(j-3)}^{2k-2-(j-1)}+...,\
 \ \ j<2m,
$$
$$ \#P^{2m}_{j}=\#T^{j-1}_{2m-2}\cdot\#T^{2k-2-(j-1)}_0+\#T^{j-1}_{(2m-2)-2}\cdot\#T^{2k-2-(j-1)}_2+...,\
 \ \ 2m\le j,\ j\ is\ odd,
$$
\[
 \lb{l222} \#P^{2m}_{j}=\#T^{j-1}_{(2m-2)-1}\cdot\#T^{2k-2-(j-1)}_1+\#T^{j-1}_{(2m-2)-3}\cdot\#T^{2k-2-(j-1)}_3+...,\
 \ \ 2m\le j,\ j\ is\ even.
\]
Using \er{2f}, \er{Pev}, \er{Pod}, we have
\[
 \lb{l2222}
 \vp_{n+1}^0(\l)=\vp_{n+1}(\l,0)=\l^n\#T_n^n-\l^{n-2}\#T^n_{n-2}+...
\]
Then \er{9999}, \er{l2222} and \er{l222}   yield
$$
 A_{1,j}=2^{-\d_{k,j}}\#T^{j-1}_{j-1}\cdot\#T^{2k-2-(j-1)}_{2k-2-(j-1)}=
 2^{-\d_{k,j}}\#P^{2k}_j=a_{1,j},
$$
$$
 -A_{2,j}=2^{-\d_{k,j}}\left(\#T^{j-1}_{j-1}\cdot\#T^{2k-2-(j-1)}_{2k-2-(j-1)-2}+
 \#T^{j-1}_{j-1-2}\cdot\#T^{2k-2-(j-1)}_{2k-2-(j-1)}\right)=
 2^{-\d_{k,j}}\#P^{2k-2}_j=a_{2,j}
$$
and so on. \BBox

\no {\bf Acknowledgements.}
The first author
E.Korotyaev was partly supported by DFG project BR691/23-1.
The authors would like also to thank Markus Klein for useful discussions.

\noindent {\bf References}

\no [AS] M. Abramowitz\ and\ I. A. Stegun,  Handbook of
mathematical functions with formulas, graphs, and mathematical
tables, U.S. Government Printing Office, Washington, D.C., 1964.

\no  [BGGK]  B\"attig, D.; Grebert, B.; Guillot, J.-C.; Kappeler,
T. Fibration of the phase space of the periodic Toda
   lattice. J. Math. Pures Appl. (9) 72 (1993), no. 6, 553--565.

\no [GT] Garnett J., Trubowitz E.: Gaps and bands of one
dimensional
 periodic Schr\"odinger operators. Comment. Math. Helv. 59, 258-312 (1984)

\no [GKT] Gieseker, D.; Kn\"orrer, H.; Trubowitz, E. The geometry
of algebraic Fermi curves. Perspectives in
   Mathematics, 14. Academic Press, Inc., Boston, MA, 1993.

\no [KK] Kargaev P., Korotyaev E.: Inverse Problem for the Hill
Operator, the Direct Approach.  Invent. Math., 129(1997), no. 3,
567-593

\no [K1] Korotyaev, E.: The inverse problem for the Hill operator.
I. Internat. Math. Res. Notices 3(1997), 113--125.

\no [K2] Korotyaev, E.: The inverse problem and trace formula for
the Hill operator, II.  Math. Z. 231(1999), 345-368

\no [K3] Korotyaev, E.: Characterization of the spectrum of
Schr\"odinger operators with periodic distributions, Int. Math.
Res. Not. 37(2003), 2019-2031

\no [K4] Korotyaev, E.: Gap-length mapping for periodic Jacobi
matrices, preprint 2004

\no [KKr]  Korotyaev, E., Krasovsky, I.: Spectral estimates for
periodic Jacobi matrices. Comm. Math. Phys. 234 (2003), no. 3,
517--532.

\no [KKu1] Korotyaev, E., Kutsenko, A.: Inverse problem for the
periodic Jacoby matrices, preprint 2004

\no [KKu2] Korotyaev, E., Kutsenko, A.: Inverse problem for the
discrete 1D Schr\"odinger operator with large periodic potential,
preprint 2004

\no [L] Last, Y. On the measure of gaps and spectra for discrete
$1$D Schrodinger operators. Comm. Math. Phys. 149 (1992), no. 2,
347--360.

\no [M]  Marchenko V.: Sturm-Liouvill operator and applications.
Basel: Birkhauser 1986.

\no  [MO1] V. Marchenko, I. Ostrovski:
A characterization of the spectrum
 of the Hill operator. Math. USSR  Sbornik  26, 493-554 (1975).

\no  [MO2]  V. Marchenko, I. Ostrovski : Approximation of periodic
by finite-zone potentials. Selecta Math. Sovietica. 1987, 6, No 2, 101-136.

\no [Pe] Perkolab, L.: An inverse problem for a periodic Jacobi
matrix. (Russian) Teor. Funktsii Funktsional. Anal. i
   Prilozhen. 42(1984), 107-121

\no [Te] Teschl, G.: Jacobi operators and completely integrable
nonlinear lattices. Mathematical Surveys
   and Monographs, 72. American Mathematical Society,
   Providence, RI, 2000

\no [vM] van Moerbeke, P.: The spectrum of Jacobi matrices.
Invent. Math. 37 (1976), no. 1, 45--81

\no [vMou] van Mouche, P. Spectral asymptotics of periodic discrete Schrodinger operators. I. Asymptotic Anal. 11 (1995), no. 3, 263--287

\end{document}